\input amstex 
\documentstyle{amsppt}
\input bull-ppt
\keyedby{bull289/lic}
\define\osl{\operatorname{SL}}
\define\dv{\operatorname{div}}
\define\grad{\operatorname{grad}}

\topmatter
\cvol{27}
\cvolyear{1992}
\cmonth{July}
\cyear{1992}
\cvolno{1}
\cpgs{134-138}
\ratitle
\title One Cannot Hear the Shape of a Drum \endtitle
\author Carolyn Gordon, David L. Webb, and Scott 
Wolpert\endauthor
\shortauthor{Carolyn Gordon, D. L. Webb, and Scott Wolpert}
\shorttitle{One Cannot Hear the Shape of a Drum}
\address \RM{(C. Gordon and D. L. Webb)} Department of 
Mathematics,
Dartmouth College, Hanover, New Hampshire 03755 \endaddress
\cu \RM{C. Gordon and D. L. Webb:} Department of 
Mathematics,
Washington University, St.~Louis, Missouri 63130\endcu
\ml C. Gordon: carolyn.gordon\@dartmouth.edu\newline
\indent{\it E-mail address\/}:  {\rm D. L. Webb: 
david.webb\@dartmouth.edu}\endml
\address \RM{(S. Wolpert)} Department of Mathematics, 
University of
Maryland, College Park, Maryland 20742\endaddress
\ml S. Wolpert: saw\@anna.umd.edu\endml
\date July 11, 1991 and, in revised form, November 5, 
1991\enddate
\thanks The authors gratefully acknowledge partial support 
from NSF grants\endthanks
\abstract We use an extension of Sunada's theorem to 
construct a
nonisometric pair of isospectral simply connected domains 
in the
Euclidean plane, thus answering negatively Kac's question, 
``can
one hear the shape of a drum?'' In order to construct 
simply 
connected examples, we exploit the observation that an 
orbifold 
whose underlying space is a simply connected manifold with 
boundary 
need not be simply connected as an orbifold.\endabstract
\endtopmatter

\document
\heading 1. Kac's question\endheading
Let $(M,g)$ be a compact Riemannian manifold with 
boundary. Then $M$ has
a Laplace operator $\Delta$, defined by 
$\Delta(f\,)=-\dv(\grad f\,)$,
that acts on smooth functions on $M$. The {\it spectrum\/} 
of $M$ is 
the sequence of eigenvalues of $\Delta$. Two Riemannian 
manifolds are
{\it isospectral\/} if their spectra coincide (counting 
multiplicities).
A natural question concerning the interplay of analysis 
and geometry is:
must two isospectral Riemannian manifolds actually be 
isometric? (When
$M$ has nonempty boundary, one can consider the {\it 
Dirichlet spectrum},
i.e., the spectrum of $\Delta$ acting on smooth functions 
that vanish on
the boundary, or the {\it Neumann spectrum}, that of 
$\Delta$ acting on
functions with vanishing normal derivative at the 
boundary.) If $M$ is a
domain in the Euclidean plane then the Dirichlet 
eigenvalues of $\Delta$ are
essentially the frequencies produced by a drumhead shaped 
like
$M$, so the question has been phrased by Bers and Kac 
\cite{16} (the latter
attributes the problem to Bochner) as ``can one hear the 
shape of a drum?''
We answer this question negatively by constructing a pair 
of nonisometric
simply connected plane domains that have both the same 
Dirichlet spectra
and the same Neumann spectra. The domains are depicted in 
Figure 1. 

The simple idea exploited here also permits us to 
construct the following:
(1) a pair of isospectral flat surfaces (with boundary) 
one of which has a
unit-length closed geodesic while the other has only a 
unit-length closed
billiard trajectory; (2) a pair of isospectral potentials 
for the Schr\"odinger
operator on

\fgh{8.5 pc}\caption{\smc Figure 1}

\noindent a plane domain (using a technique of Brooks 
\cite5); (3) a pair
of isospectral, nonisometric domains in the hyperbolic 
plane; and (4) a
pair of isospectral, nonisometric domains in the 2-sphere.

Weyl's proof \cite{23} that the area of a plane domain is 
determined
by the spectrum led to speculation that perhaps the shape 
of a plane
domain (or more generally, of a Riemannian manifold) is 
audible. The
latter was refuted by Milnor \cite{17}, who exhibited a 
pair of isospectral,
nonisometric 16-dimensional tori. Other examples followed,
including (among others) isospectral pairs of Riemann 
surfaces constructed
by Vign\'eras \cite{22}, Buser \cite{7, 8}, Brooks \cite4, 
and Brooks-Tse
\cite6; pairs of lens spaces produced by Ikeda \cite{15}; 
pairs of domains
in $\bold R^4$, due to Urakawa \cite{21}; and continuous 
families of 
isospectral metrics on solvmanifolds constructed by 
Gordon-Wilson
\cite{14} and DeTurck-Gordon \cite{11, 12}. However, Kac's 
question
concerning plane domains has remained open.

As will be clear from the discussion of Sunada's theorem 
below, most
known pairs of isospectral manifolds have a common 
Riemannian cover.
Thus it is also of interest to exhibit simply connected 
isospectral
manifolds.

\heading 2. Sunada's Theorem\endheading
Although the early examples of isospectral manifolds 
seemed rather 
{\it ad hoc}, a coherent explanation for most of them has 
since been
provided. Sunada \cite{19} introduced a general method for 
constructing
pairs of isospectral manifolds with a common finite 
covering:

\thm\nofrills{Theorem \RM{(Sunada).}}\ Let $M$ be a 
Riemannian manifold
upon which a finite group $G$ acts by isometries\RM; let 
$H$ and $K$ be
subgroups of $G$ that act freely. Suppose that $H$ and $K$ 
are {\it almost
conjugate}, i.e., there is a bijection $f\: H\to K$ 
carrying every
element  $h$ of $H$ to an element $f(h)$ of $K$ that is 
conjugate
in $G$ to $h$. Then the quotient manifolds 
$M_1=H\backslash M$ and 
$M_2=K\backslash M$ are isospectral.
\ethm

Choosing conjugate subgroups in the above theorem yields 
isometric manifolds,
so one seeks a finite group with a pair of almost 
conjugate but nonconjugate
subgroups. The algebraic condition can be restated as: the 
representations 
$L^2(H\backslash G)$ and $L^2(K\backslash G)$ are 
unitarily equivalent,
although $H\backslash G$ and $K\backslash G$ are 
inequivalent as $G$\<-sets.

B\'erard \cite1 gave a new proof of Sunada's theorem, by 
noting the following:

\thm{Proposition} Let $G$ be a group with subgroups $H$ 
and $K$, and suppose
that $T\: L^2(H\backslash G)$ $\to$ $L^2(K\backslash G)$ is a 
unitary intertwining
operator. Let $V$ be a Hilbert space on which $G$ acts 
unitarily. Then $T$
induces an isometry $V^H\to V^K$. \RM(Here $V^H$ denotes 
the subspace of
$V$ consisting of the $H$\<-fixed points.\RM)\ethm

Note that if $G$ acts by isometries on $M$ and hence on 
$V=L^2(M)$, then
$V^H=L^2(H\backslash M)$ since the $H$\<-invariant 
functions on $M$ are
precisely those functions that descend to the quotient 
manifold, and likewise
$V^K=L^2(K\backslash M)$. B\'erard then used the 
proposition along with a
variational characterization of eigenvalues of the 
Laplacian to recover
Sunada's theorem. He also pointed out that the assumption 
that $H$ and
$K$ act freely is not necessary. In this case the 
quotients $H\backslash
M$ and $K\backslash M$ are orbifolds, as discussed below, 
although their
underlying spaces may be manifolds with boundary.

\heading 3. Orbifolds\endheading
An {\it orbifold\/} is a space locally modelled on the 
orbit space of a
finite group acting on $\bold R^n$; for a precise 
definition, see \cite{20}
or \cite{18}. In particular, the quotient space 
$O=G\backslash M$ of a 
manifold $M$ by a group $G$ acting properly 
discontinuously is an orbifold.
If $G$ acts freely, then $O$ is a manifold; otherwise $O$ 
may have a singular
set arising from fixed points of the action of $G$. 

There are modified definitions of the fundamental group of 
an orbifold and
of orbifold covering maps. The important feature for our 
purposes is
that the underlying space $|O|$ of an orbifold $O$ may 
have no ordinary
proper covering spaces, although $O$ has proper coverings 
in the orbifold
sense; thus an orbifold with a simply connected underlying 
space need not
be simply connected as an orbifold. 

Viewing $G$\<-invariant functions on $M$ as functions on 
$G\backslash M$,
one defines the spectrum of a quotient orbifold 
$G\backslash M$ as
the eigenvalue spectrum of the Laplace  operator acting on 
the space
$L^2(M)^G$ of $G$\<-invariant functions on $M$. In 
particular, if $G$
is a group that acts by isometries on a Riemannian 
manifold $M$ and if
$H$, $K$ are almost conjugate subgroups of $G$, perhaps 
acting with 
fixed points, then from the above discussion one obtains 
isospectral
orbifolds $O_1=H\backslash M$ and $O_{\,2}=K\backslash M$. 

\heading 4. Construction of isospectral simply connected 
manifolds\endheading
We now utilize the above observations to construct 
isospectral simply
connected plane domains. We use the discussion of \S2 to 
produce an 
isospectral pair $O_1$, $O_{\,2}$ of 2-orbifolds with 
boundary by modifying
a construction due to Buser \cite8 of an isospectral pair 
$M_1$, $M_2$
of flat surfaces with boundary. Buser's surfaces are 
constructed as covers
of a bordered surface $M_0$ using a pair of almost 
conjugate subgroups of
$G=\osl_3(\bold F_2)$ and a representation of $\pi_1(M_0)$ 
in $G$, and our
orbifolds are similarly constructed as covers of an 
orbifold $O_0$, using
the orbifold notion of fundamental group, the 
corresponding theory of
orbifold coverings, and a representation of $\pi_1(O_0)$ 
in a split
extension of $G$ by $\bold Z/2\bold Z$; indeed, the 
orbifold $O_i$ is
the quotient by an involutive isometry of Buser's manifold 
$M_i$,
$i=0,1,2$. We observe that the Neumann orbifold spectrum 
of $O_i$ is
precisely the Neumann spectrum of the underlying manifold 
$|O_i|$;
thus the underlying spaces $|O_1|$ and $|O_{\,2}|$ are 
Neumann-isospectral
manifolds with boundary. $O_1$ and $O_{\,2}$ have a common 
cover in the
orbifold sense, but not in the usual sense; this common 
cover $O$ is
the quotient by an involutive isometry of a common cover 
$M$ of Buser's
surfaces $M_1$ and $M_2$. The underlying spaces $|O_1|$ 
and $|O_{\,2}|$
are simply connected plane domains. We deduce the 
Dirichlet isospectrality
of $|O_1|$ and $|O_{\,2}|$ by exploiting the Dirichlet 
isospectrality of the
double covers $M_1$ and $M_2$, which corresponds to 
isospectrality of
the plane domains with mixed boundary conditions: 
Dirichlet conditions
on the orbifold boundary and Neumann conditions on the 
singular set; this
observation, together with the obvious decomposition of 
any eigenfunction
into an involution--invariant eigenfunction and an 
involution--anti-invariant
eigenfunction establishes the Dirichlet isospectrality. 
Details will
appear elsewhere.

\heading Acknowledgment\endheading 
The first two authors wish to thank the Institute Fourier 
for its
hospitality during the period when some of this work was 
done.
In particular, they wish to thank Pierre B\'erard, Gerard 
Besson,
Bob Brooks, and Yves Colin de Verdi\`ere for stimulating 
discussions.

\enddocument